\documentclass[12pt, reqno]{amsart}
\usepackage{amsmath, amsthm, amscd, amsfonts, amssymb}
\usepackage[mathscr]{eucal}
\usepackage[bookmarksnumbered, colorlinks, plainpages]{hyperref}
\hypersetup{colorlinks=true,linkcolor=red, anchorcolor=green, citecolor=cyan, urlcolor=red, filecolor=magenta, pdftoolbar=true}

\usepackage[pdftex]{graphicx}

\newtheorem{theorem}{Theorem}

\newtheorem{proposition}[theorem]{Proposition}

\theoremstyle{remark}
\newtheorem{remark}{Remark}

\DeclareMathOperator{\sinch}{sinch}

\newcommand{\m}{\mathfrak m}

\begin{document}
\setcounter{page}{1}

\title[Calculating function of a matrix]{Calculating a function of a matrix\\ with a real spectrum}

\author{P.~Kubel\'{i}k}
 \address{J.~Heyrovsk\'{y} Institute of Physical Chemistry, Academy of Sciences of the Czech Republic, Dolej\v{s}kova 3, 18223 Prague 8, Czech Republic}
\email{\textcolor[rgb]{0.00,0.00,0.84}{petr.kubelik@jh-inst.cas.cz}}

\author{V.~G. Kurbatov}
 \address{Department of System Analysis and Control,
Voronezh State University\\ 1, Universi\-tet\-skaya Square, Voronezh 394018, Russia}
\email{\textcolor[rgb]{0.00,0.00,0.84}{kv51@inbox.ru}}

\author{I.~V. Kurbatova}
 \address{Department of Software Development and Information Systems Administration,
Vo\-ro\-nezh State University\\ 1, Universitetskaya Square, Voronezh 394018, Russia}
\email{\textcolor[rgb]{0.00,0.00,0.84}{irakurbatova@gmail.com}}

\subjclass{Primary 65F60; Secondary 97N50}

\keywords{matrix function, interpolating polynomial, polynomial interpolation, Chebyshev nodes, divided differences, reordering, significant decimal digits}

\date{\today}

\begin{abstract}
Let $T$ be a square matrix with a real spectrum, and let $f$ be an analytic function. The problem of approximate calculation of $f(T)$ is discussed. Applying the Schur triangular decomposition and the reordering, one can assume that $T$ is triangular and its diagonal entries $t_{ii}$ are arranged in increasing order. To avoid calculations using the differences $t_{ii}-t_{jj}$ with close (including equal) $t_{ii}$ and $t_{jj}$, it is proposed to represent $T$ in a block form and calculate the two main block diagonals using interpolating polynomials. The rest of the entries of $f(T)$ can be calculated by means of the Parlett recurrence algorithm. It is also proposed to perform scalar operations (such as building of interpolating polynomials) with an enlarged number of decimal digits.
\end{abstract}

\maketitle

\section{Introduction}\label{s:Introduction}
The problem~\cite{Godunov94:ODE:eng,Higham08} of calculating analytic functions $f$ of a matrix $R$ is of constant interest since it has at least one extensive application: the solution of a linear differential equation $\dot x=Rx$ with a matrix coefficient $R$ can be expressed in terms of the function $e^{Rt}$. We only mention some other applications: the cosine~\cite{Grimm-Hochbruck08,Hargreaves-Higham05,Sastre-Ibane-Alonso-Peinado-Defez17}, the power function~\cite{Burrage-Hale-Kay12,Higham-Lin2011}, the sign function~\cite{Bai-Demmel98,Higham08,Kenney-Laub95}, Green's function~\cite{Godunov94:ODE:eng,Kurbatov-Kurbatova-CMAM18,Kurbatova-Pechkurov20}, the gamma and beta functions~\cite{Cardoso-Sadeghi18,Jodar-Cortes,Schmelzer-Trefethen07}, the Lambert function~\cite{Jarlebring07}. Usually, a function of a matrix can be calculated only approximately. The main objectives are (i) consideration of matrices of large size, (ii) obtaining the final result with high accuracy, and (iii) fast computations. In this paper, we propose a new idea which expands the possibilities of accurate calculations.

We present a numerical algorithm for the approximate calculation of an analytic function $f$ of a square matrix $R$ with a real spectrum.
A trivial (non-interesting from the point of view of this paper) example of a matrix with a real spectrum is a Hermitian matrix. A more substantive example of a matrix with a real spectrum is~\cite[12.69]{Voevodin-Kuznetsov:rus-eng} a matrix of the form $F^{-1}G$, where $F$ and $G$ are Hermitian and $F$ is positive definite. The matrix $F^{-1}G$ arises when one transforms the linear pencil~\cite{Ikramov91a:eng} $\lambda\mapsto\lambda F-G$ to the form $\lambda\mapsto\lambda\mathbf1-F^{-1}G$, where $\mathbf1$ is the identity matrix. In its turn, linear pencils $\lambda\mapsto\lambda F-G$ with positive definite matrices $F$ and $G$ describe~\cite{Seshu-Reed61:eng,Vlach-Singhal94,Kurbatov-Oreshina04} $RC$-circuits and some other passive systems.

First, we construct the Schur decomposition, i.~e. the representation $R=QTQ^H$, where $Q$ is unitary and $T$ is triangular.
We recall that the Schur decomposition is computed with high backward stability by the
$QR$-algorithm~\cite{Golub-Van_Loan13:eng,Householder64,Voevodin-Kuznetsov:rus-eng,Wilkinson88}. We also note that the Schur decomposition is used as a beginning step in many numerical spectral algorithms. Clearly, one has
\begin{equation*}
f(R)=Qf(T)Q^H.
\end{equation*}
Thus, the problem of calculating $f(R)$ is reduced to calculating $f(T)$ for a triangular matrix $T$. (We recall that the diagonal entries of $T$ are the eigenvalues of both $T$ and $R$.)
Then we apply an additional unitary similarity transformation that arranges the diagonal entries of $T$ in the increasing order (see Section~\ref{s:reodering}).

The most known and widely used algorithm for calculating an analytic function $f$ of an upper triangular matrix $T$ is due to Parlett~\cite{Parlett74,Parlett76}, see also its expositions in~\cite[9.1.4]{Golub-Van_Loan13:eng} and~\cite[p.~85]{Higham08}. If the matrix $T$ has no multiple eigenvalues, the calculation of the entries $f_{ij}$ of the matrix $F=f(T)$ (where $T$ is upper triangular) can be performed by the formulae
\begin{equation}\label{e:Parlett:single}
\begin{split}
f_{ii}&=f(t_{ii}),\\
f_{ij}&=t_{ij}\frac{f_{jj}-f_{ii}}{t_{jj}-t_{ii}}
+\sum_{k=i+1}^{j-1}\frac{t_{ik}f_{kj}-f_{ik}t_{kj}}{t_{jj}-t_{ii}},\qquad i<j.
\end{split}
\end{equation}
It is seen from~\eqref{e:Parlett:single} that the diagonals of the matrix $F$ can be calculated one after another. The usage of the Schur decomposition and then application of formulae~\eqref{e:Parlett:single} is usually called~\cite{Golub-Van_Loan13:eng,Higham08} the \emph{Schur--Parlett algorithm}.

Unfortunately, formulae~\eqref{e:Parlett:single} do not work if $T$ has multiple eigenvalues, i.~e. $t_{ii}=t_{jj}$ for some $i\neq j$. A solution to this problem was also proposed by Parlett~\cite{Parlett74}, see also~\cite{McCurdy-Ng-Parlett84,Golub-Van_Loan13:eng,Higham08}: the diagonal entries should be reordered so that the matrix can be represented in a block form
with equal eigenvalues $t_{ii}$ included in common blocks; after that, a block analogue~\cite{Golub-Van_Loan13:eng,Higham08} of formulae~\eqref{e:Parlett:single} can be applied.

A similar problem arises when $T$ has very close (though unequal) diagonal entries $t_{ii}$ and $t_{jj}$. In this case, the calculation of the denominators $t_{jj}-t_{ii}$ in~\eqref{e:Parlett:single} causes a severe loss of accuracy. For a detailed discussion of this phenomenon, see~\cite{Moler-Van_Loan78,Moler-Van_Loan03}. A partial solution of this problem was proposed in~\cite{Davies-Higham03}. Close eigenvalues should be united into clusters. After that the diagonal entries of the triangular matrix should be reordered in such a way that the eigenvalues from the same cluster follow one after another. Thus, the matrix is divided into blocks so that the eigenvalues from the same cluster fall into the same block. Finally, the function $f$ of diagonal blocks (having small spectrum) can be calculated using the Taylor (or interpolating) polynomial; then the remaining blocks of $F$ can be calculated using the block version of~\eqref{e:Parlett:single}; note that the absence of $t_{jj}-t_{ii}$ with close $t_{jj}$ and $t_{ii}$ in the neighboring blocks requires that the clusters be separated from each other. A modification of this algorithm that does not use a reordering was discussed in~\cite{Kurbatov-Kurbatova-LMA16}.

In this paper, we consider the case when partitioning into well-separated clusters
is impossible: the spectrum of $R$ uniformly (in the statistical sense) fills a segment of the real line. The main idea is as follows. Since we assume that the spectrum of $R$ is real, we can rearrange the diagonal entries of $R$ in increasing order. We again divide the spectrum of $R$ into clusters consisting of close numbers. Then we unite each pair of neighboring clusters and consider the diagonal blocks that correspond to these pairs (thus each cluster except the two extreme ones is involved into two diagonal blocks), see Section~\ref{s:partitioning} for details. For example, let the spectrum of $R$ be contained in $[0,3]$ and we form 3 cluster: eigenvalues that are contained in $[0,1]$, eigenvalues that are contained in $[1,2]$, and eigenvalues that are contained in $[2,3]$; then we consider two (overlapping) diagonal blocks that correspond to the eigenvalues from $[0,2]$ and $[1,3]$. We calculate the function $f$ of diagonal blocks (note that the diagonal blocks have relatively small spectra) by means of the Taylor polynomials or another special tool. Unfortunately, since the diagonal blocks intersect, their common entries are calculated twice. Actually, we use a modification of the described idea that allows one to avoid such duplication. We calculate the rest of entries of $F$ using~\eqref{e:Parlett:single}; since the rest of the entries are located far from the main diagonal, the differences $t_{jj}-t_{ii}$ with close $t_{jj}$ and $t_{ii}$ do not occur.

To calculate the action of $f$ on blocks with a small (real) spectrum, we use the Newton interpolating polynomials with the Chebyshev nodes taken as the points of interpolation. The construction of the Newton polynomial requires the calculation of divided differences corresponding to close points of interpolation. This operation results in a loss of accuracy. We propose performing it with an enlarged number of decimal digits. No significant time is spent on such calculations (in comparison with matrix operations). By the way, the experiment shows that the longest operation is the calculation according to formula~\eqref{e:Parlett:single} of the entries that are far from the main diagonal.

The paper is organized as follows. In Section~\ref{s:reodering}, we recall the Schur decomposition and the reordering, and introduce some notation. In Section~\ref{s:partitioning}, we begin the description of our algorithm, represent the matrix in the block form, and describe the simplified version of the algorithm. In Section~\ref{s:small spectrum}, we recall and discuss methods of calculating blocks from the main block diagonal. In Section~\ref{s:Paterson-Stockmeyer}, we recall the Paterson--Stockmeyer algorithm that accelerates the substitution of a matrix into a polynomial. Section~\ref{s:DD} describes an algorithm for calculating the second block diagonal. In Section~\ref{s:algorithm}, we list the sequence of execution of individual parts of our algorithm. In Section~\ref{s:numerical example}, numerical experiments are presented.

\section{Reordering and other preliminaries}\label{s:reodering}
In this section, we describe the preliminary stage of our algorithm.

Let $N$ be a natural number and $R$ be a complex matrix of the size $N\times N$ with a real spectrum, and $f$ be a function analytic in a neighborhood of the spectrum of $R$. In this paper, we present an algorithm of the approximate calculation of the matrix $f(R)$. Our algorithm is interesting in the case when the spectrum $\sigma(R)$ of $R$ fills (in a statistical sense) a segment of the real line $\mathbb R$ without essential gaps. The algorithm can be modified to the case when the spectrum $\sigma(R)$ lies on a continuous curve in $\mathbb C$ without self-intersections.

First, we compute the Schur decomposition $R=QTQ^H$, where the matrix $Q$ is unitary and $T$ is triangular.
We recall that the Schur decomposition is performed with high backward stability by the
$QR$-algorithm~\cite{Golub-Van_Loan13:eng,Householder64,Voevodin-Kuznetsov:rus-eng,Wilkinson88}.

Our algorithm needs that the diagonal entries of $T$ (the eigenvalues of both $T$ and $R$) be arranged in order.
Usually, the $QR$-algorithm arranges the diagonal entries of $T$ in an almost descending order of modules. For our notation, it is convenient to have the increasing order. Therefore, we apply the $QR$-algorithm to the matrix $R_1=-R+\alpha\mathbf1$, where $\mathbf1$ is the identity matrix and $\alpha\in\mathbb R$ is greater than the largest eigenvalue of $R$. For example, one may take for $\alpha$ a number that is greater than the norm $\Vert R\Vert_{1\to1}=\max_{j}\sum_{i=1}^m|t_{ij}|$ of $R$. We calculate the decomposition $R_1=QT_1Q^H$, where $Q$ is unitary and $T_1$ is triangular. Then we set $T=T_1+\alpha\mathbf1$ and arrive at the Schur decomposition $R=QTQ^H$ of the matrix $R$, in which the diagonal entries of $T$ are arranged in almost increasing order.

Since the order of the diagonal entries $t_{ii}$ may be not fully increasing, we correct it using the procedure of \emph{reordering} (see, e.g.,~\cite{Bai-Demmel93}, \cite[p.~49]{Ikramov84:rus-eng}), i.~e. we apply several unitary similarity transformations $T'=UTU^H$, where $U$ is unitary, that permutes two adjacent diagonal entries $t_{i-1,i-1}$ and $t_{i,i}$ but preserves the triangular structure of $T$. For completeness, we present the description of the unitary similarity transformation $U$: the entries of $U$ coincide with the entries of the identity matrix except for the block
\begin{equation*}
\begin{pmatrix}
u_{i-1,i-1}&u_{i-1,i}\\
u_{i,i-1}&u_{i,i}\\
\end{pmatrix}=
\begin{pmatrix}
e^{-i\alpha}\cos\vartheta&-\sin\vartheta\\
\sin\vartheta&e^{i\alpha}\cos\vartheta\\
\end{pmatrix},
\end{equation*}
where the parameters $\alpha$ and $\vartheta$ are real. The inverse (conjugate transpose) of this block has the form
\begin{equation*}
\begin{pmatrix}
e^{-i\alpha}\cos\vartheta&-\sin\vartheta\\
\sin\vartheta&e^{i\alpha}\cos\vartheta\\
\end{pmatrix}^{-1}=
\begin{pmatrix}
e^{i\alpha}\cos\vartheta&\sin\vartheta\\
-\sin\vartheta&e^{-i\alpha}\cos\vartheta\\
\end{pmatrix}.
\end{equation*}
The first admissible pair of the parameters $\alpha$ and $\vartheta$ are defined from the formulae
\begin{align*}
\alpha_1&=\arg t_{i-1,i}-\arg(t_{i-1,i-1}-t_{i,i}),\\
\cos\vartheta_1&=\frac{|t_{i-1,i}|}{\sqrt{|t_{i-1,i-1}-t_{i,i}|^2+|t_{i-1,i}|^2}},\\
\sin\vartheta_1&=\frac{|t_{i-1,i-1}-t_{i,i}|}{\sqrt{|t_{i-1,i-1}-t_{i,i}|^2+|t_{i-1,i}|^2}}.
\end{align*}
The second admissible pair is
\begin{equation*}
\alpha_2=\alpha_1+\pi,\qquad \vartheta_2=-\vartheta_1.
\end{equation*}
The two possible values of $\alpha$ and $\vartheta$ lead only to different signs in the entries of $T'=UTU^H$. Therefore, it does not matter which pair to use.

\section{Partitioning}\label{s:partitioning}
Below we assume that the diagonal entries $t_{ii}$ of the triangular matrix $T$ are real and arranged in ascending order.

Let $a$ be the minimum of the diagonal entries $t_{ii}$, and let $b$ be greater than the maximum of the diagonal entries $t_{ii}$. Thus, the spectrum of $T$ is contained in the half-open interval\footnote{The usage of a half-open interval is due only to the uniformity of the notation. In practice, it is more convenient to use a closed interval $[a,b]$.} $[a,b)$.

We choose the parameter $\rho>0$. The algorithm depends on this parameter. In our numerical examples in Section~\ref{s:numerical example} $\rho=2$. The meaning of the parameter $\rho$ is as follows: we consider the loss of accuracy to be acceptable if formula~\eqref{e:Parlett:single} is applied only in the case $|t_{jj}-t_{ii}|\ge\rho$.

We divide $[a,b)$ into parts of the length $\rho$:
\begin{equation*}
[a,b)\subseteq[a,a+\rho)\cup[a+\rho,a+2\rho)
\cup\dots\cup[a+(n-1)\rho,a+n\rho).
\end{equation*}
Thereby, we assume that $b\in[a+(n-1)\rho,a+n\rho)$. We divide the set $t_{11}$,
$t_{22}$, \dots, $t_{NN}$ of diagonal entries of $T$ into $n$ clusters:
\begin{itemize}
 \item the first cluster consists of diagonal entries such that
 \begin{equation*}
 t_{11},t_{22},\dots,t_{k_1k_1}\in[a,a+\rho);
\end{equation*}
 \item the second cluster consists of diagonal entries such that
 \begin{equation*}
 t_{k_1+1,k_1+1},t_{k_1+2,k_1+2},\dots,t_{k_2k_2}\in[a+\rho,a+2\rho);
 \end{equation*}
 \item \dots;
 \item the $n$-th cluster consists of diagonal entries such that
 \begin{equation*}
 t_{k_{n-1}+1,k_{n-1}+1},t_{k_{n-1}+2,k_{n-1}+2},\dots,t_{k_n k_n}\in[a+(n-1)\rho,a+n\rho).
 \end{equation*}
\end{itemize}
We divide the matrix $T$ into blocks according to this clusterization:
\begin{equation}\label{e:block form of T}
T=\begin{pmatrix}
T_{1,1} & T_{1,2} & \dots & T_{1,n-1} & T_{1,n} \\
0 & T_{2,2} & \dots & T_{2,n-1} & T_{2,n} \\
\vdots & \vdots & \ddots & \vdots & \vdots \\
0 & 0 & \dots & T_{n-1,n-1} &  T_{n-1,n} \\
0 & 0 & \dots & 0 & T_{n,n}
 \end{pmatrix}.
\end{equation}

Now we come to the main idea of the paper. To begin with, we describe its simplest realization. We fill the first two diagonals of the block matrix $F$ by applying the function $f$ to the double blocks
\begin{equation}\label{e:double blocks}
\begin{pmatrix}
T_{1,1} & T_{1,2}\\
0 & T_{2,2}\end{pmatrix},\qquad
\begin{pmatrix}
T_{2,2} & T_{2,3}\\
0 & T_{3,3}\end{pmatrix},\qquad\dots,\qquad
\begin{pmatrix}
T_{n-1,n-1} &  T_{n-1,n} \\
0 & T_{n,n}\end{pmatrix}.
\end{equation}
Since these blocks have relatively narrow spectra, the result of the application of the function $f$ to the double blocks can be approximately calculated, e.g., by substituting them into the Taylor polynomial of $f$.
(Unfortunately, the diagonal blocks $F_{2,2}$, \dots, $F_{n-1,n-1}$ of $F=f(T)$ are calculated twice. Of course, this leads to a loss of time. Below we eliminate this duplication.) Thus, we calculate two main diagonals of the block matrix $f(T)$. After that, the remaining entries of $F$ are calculated by formula~\eqref{e:Parlett:single}.
This is the main idea of the paper. In the rest of the paper, we discuss its improvements that work faster and more accurate. Namely, we describe how to calculate the two main block diagonals of $F$ separately.

The idea of calculating first the entries of two main block diagonals and then the other entries was used earlier in~\cite{McCurdy-Ng-Parlett84}, but for a somewhat different purpose, namely, for applying $f$ to the two-diagonal Opitz matrix, which allows one to calculate the divided differences of $f$ in order to then substitute them into the Newton interpolating polynomial for calculating $f(R)$. We use an analogous idea to calculate $f(R)$ directly.

\section{Calculation of $F_{kk}$}\label{s:small spectrum}
In this section, we discuss the ways of calculating the function $f$ of the diagonal blocks of~\eqref{e:block form of T}. By construction, the spectrum of any block $T_{kk}$ is real and contained in a segment of length $\rho$. This circumstance simplifies the problem.

We recall the following well-known statement.
  \begin{proposition}[{\rm \cite[Theorem 1.13]{Higham08}}]\label{p:inv by det}
Let a block matrix $T$ have form~\eqref{e:block form of T} and the matrix $F=f(T)$ be presented in the same block form{\rm:}
\begin{equation}\label{e:block form of F}
F=\begin{pmatrix}
F_{1,1} & F_{1,2} & \dots & F_{1,n-1} & F_{1,n} \\
0 & F_{2,2} & \dots & F_{2,n-1} & F_{2,n} \\
\vdots & \vdots & \ddots & \vdots & \vdots \\
0 & 0 & \dots & F_{n-1,n-1} &  F_{n-1,n} \\
0 & 0 & \dots & 0 & F_{n,n}
 \end{pmatrix}.
\end{equation}
Then
\begin{equation*}
F_{kk}=f(T_{kk}),\qquad k=1,2,\dots,n.
\end{equation*}
 \end{proposition}

The purpose of this section is to discuss the approximate calculation of $f(T_{kk})$.

We fix $k=1,2,\dots,n$.
Let $[\alpha_k,\beta_k]\subseteq\mathbb R$ be a segment that contains the spectrum of the block $T_{kk}$ (e.g., $\alpha_k$ is the minimal eigenvalue of $T_{kk}$ and $\beta_k$ is the maximal one.). Clearly, $\beta_k-\alpha_k\le\rho$.
We denote by $\gamma_k$ the center of the segment $[\alpha_k,\beta_k]$, i.~e. the point $\frac{\alpha_k+\beta_k}2$.

We denote by $A$ the matrix $T_{kk}-\gamma_k\mathbf1$, where $\mathbf1$ is the identity matrix, and we denote by $g$ the shifted function $f$:
\begin{equation*}
g(\lambda)=f(\lambda+\gamma_k).
\end{equation*}
Clearly, $f(T_{kk})=g(A)$ and the spectrum of $A$ is contained in
\begin{equation*}
[\hat\alpha,\hat\beta]=\bigl[-\delta_k,\delta_k\bigr],
\end{equation*}
where $\delta_k=\frac{\beta_k-\alpha_k}2$.

We denote by $m_k$ the degree of a polynomial $p$ we will use for the approximation of $g$ (the polynomial $p$ has $\m=m_k+1$ coefficients). In our numerical experiments in Section~\ref{s:numerical example}, we take $m_k=15$.

\subsection{The usage of the Taylor polynomial}\label{ss:Taylor polynomial}
Let $p=p_k$ be the Taylor polynomial of degree $m_k$ of the function $g$ about $\lambda=0$:
\begin{equation*}
p_k(\lambda)=\sum_{k=0}^{m_k}\frac{g^{(i)}(0)}{i!}\lambda^i
=\sum_{i=0}^{m_k}\frac{f^{(i)}(\gamma_k)}{i!}\lambda^i.
\end{equation*}
We can take $p_k(A)$ as an approximation for $f(A)$. The problem is how to choose the degree $m_k$ of the Taylor polynomial to ensure enough accuracy. Surely, if the Taylor series \begin{equation*}
g(A)=\sum_{i=0}^{\infty}\frac{f^{(i)}(\gamma_k)}{i!}A^i
\end{equation*}
converges, then one has the simplest estimate
\begin{equation*}
\lVert p_{k}(A)-g(A)\rVert=\biggl\lVert \sum_{i=m_k+1}^{\infty}\frac{f^{(i)}(\gamma_k)}{i!}A^i\biggr\rVert
\le\sum_{i=m_k+1}^{\infty}\frac{\bigl| f^{(i)}(\gamma_k)\bigr|}{i!}\lVert A\rVert^i.
\end{equation*}
In many cases the right hand side of this formula can be easily calculated or estimated. For a more precise but at the same time more complicated estimate, see~\cite{Mathias93a}. In order to accelerate the substitution of $T_{kk}$ into the polynomial $p_{k}$, one can use the Paterson--Stockmeyer algorithm~\cite{Paterson-Stockmeyer73}, see Section~\ref{s:Paterson-Stockmeyer}.

It is reasonable to approximate $g$ by the Taylor polynomial $p$ of such degree that $p$ approximates $g$ on $[\hat\alpha,\hat\beta]$ within 16 decimal digits (but the degree of $p$ is not large).
The following numerical results show which degree of $p$ is acceptable depending on the length of $[\hat\alpha,\hat\beta]$.
For the Taylor polynomial $p_m$ of degree $m$ of the exponential function one has
\begin{align*}
\max_{|\lambda|\le0.12}|e^\lambda-p_9(\lambda)|&=1.72*10^{-16},\\
\max_{|\lambda|\le0.34}|e^\lambda-p_{12}(\lambda)|&=1.34*10^{-16},\\
\max_{|\lambda|\le0.68}|e^\lambda-p_{15}(\lambda)|&=1.04*10^{-16},\\
\max_{|\lambda|\le1}|e^\lambda-p_{17}(\lambda)|&=1.65*10^{-16}.
\end{align*}

\subsection{The usage of the interpolating polynomial}\label{ss:interpolation polynomial}
We prefer to calculate $g(A)$ by means of an interpolating polynomial. Interpolating polynomials can give less error than Taylor ones; in other words, an interpolating polynomial of a smaller degree can ensure the same accuracy. Since the spectrum of $A$ is real, it is reasonable to use as the points of interpolation the zeroes of the Chebyshev polynomial of the first kind $T_{\m}(x)=\cos(\m\arccos x)$ scaled to the segment $[\hat\alpha,\hat\beta]$, which contains the spectrum of $A$, i.~e. the points
\begin{equation}\label{e:Chebyshev nodes}
\lambda_{i}=\frac12(\hat\alpha+\hat\beta)+\frac12(\hat\beta-\hat\alpha)\cos\dfrac{(2i-1)\pi}{2\m},\qquad i=1,2,\dots,\m.
\end{equation}
We call points~\eqref{e:Chebyshev nodes} the \emph{Chebyshev nodes of degree $m$ on the segment} $[\hat\alpha,\hat\beta]$.
We mention that the Chebyshev nodes of degree $16$ on $[-1,1]$ are the numbers
\begin{align*}
\pm\sqrt{2-\sqrt{2+\sqrt{2+\sqrt{2}}}},\qquad
\pm\sqrt{2-\sqrt{2+\sqrt{2-\sqrt{2}}}},\\
\pm\sqrt{2-\sqrt{2-\sqrt{2-\sqrt{2}}}},\qquad
\pm\sqrt{2-\sqrt{2-\sqrt{2+\sqrt{2}}}},\\
\pm\sqrt{2+\sqrt{2-\sqrt{2+\sqrt{2}}}},\qquad
\pm\sqrt{2+\sqrt{2-\sqrt{2-\sqrt{2}}}},\\
\pm\sqrt{2+\sqrt{2+\sqrt{2-\sqrt{2}}}},\qquad
\pm\sqrt{2+\sqrt{2+\sqrt{2+\sqrt{2}}}}.
\end{align*}
Thus, they can be easily calculated with any desirable accuracy.

We recall why the choice of the Chebyshev nodes is a good one. From the Lagrange form of the interpolating polynomial
\begin{equation*}
p(\lambda)=\sum_{i=1}^{\m} g(\lambda_i)l_{i}(\lambda),
\end{equation*}
where
\begin{equation*}
l_i(\lambda)=\prod\limits_{\substack{j=1\\j\neq i}}^{\m}\frac{\lambda-\lambda_j}{\lambda_i-\lambda_j}
\end{equation*}
are the Lagrange basis polynomials, it follows that
\begin{equation}\label{e:p le Lambda f}
\max_{\lambda\in[\hat\alpha,\hat\beta]}|p(\lambda)|
\le\Lambda_{\{\lambda_i\}}\max_{\lambda\in[\hat\alpha,\hat\beta]}|g(\lambda)|,
\end{equation}
where
\begin{equation*}
\Lambda_{\{\lambda_i\}}=\max_{\lambda\in[\hat\alpha,\hat\beta]}\sum_{i=1}^{\m}|l_{i}(\lambda)|.
\end{equation*}
The number $\Lambda_{\{\lambda_i\}}$ is called~\cite[p.~212]{Babenko:eng:2nd} the \emph{Lebesgue constant} that corresponds to the segment $[\hat\alpha,\hat\beta]$ and the points of interpolation $\lambda_i$, $i=1,2,\dots,\m$. From~\eqref{e:p le Lambda f} it follows that
\begin{equation}\label{e:norm(f-p)}
\max_{\lambda\in[\hat\alpha,\hat\beta]}|g(\lambda)-p(\lambda)|
\le\max_{\lambda\in[\hat\alpha,\hat\beta]}|g(\lambda)|+\max_{\lambda\in[\hat\alpha,\hat\beta]}|p(\lambda)|
\le(1+\Lambda_{\{\lambda_i\}})\max_{\lambda\in[\hat\alpha,\hat\beta]}|g(\lambda)|.
\end{equation}

We recall that $\m$ points of interpolation generate the interpolating polynomial $p$ of degree $m=\m-1$.
We denote by $E_{m}(g)$ the error of the best uniform approximation on $[\hat\alpha,\hat\beta]$ to a function $g$ by polynomials of degree at most $m$, i.~e.
\begin{equation*}
E_{m}(g)=\min_{p\in\mathcal P_{m}}\max_{\lambda\in[\hat\alpha,\hat\beta]}|g(\lambda)-p(\lambda)|,
\end{equation*}
where the minimum is taken over the set $\mathcal P_{m}$ of all polynomials $p$ of degree at most $m$. It is known~\cite{Bernstein1914,Varga-Karpenter86} that $\lim_{m\to\infty}2m\,E_{2m}(f)\approx0.280$ for any continuous real function $f$.
The Lebesgue constant shows to what extent the interpolating polynomial $p$ of degree $m$ gives a worse approximation than the polynomial of best approximation (of the same degree), namely, the following estimate~\cite[p.~212, formula (5)]{Babenko:eng:2nd} holds:
\begin{equation*}
\max_{\lambda\in[\hat\alpha,\hat\beta]}|g(\lambda)-p(\lambda)|\le(1+\Lambda_{\{\lambda_i\}})E_{m}(g).
\end{equation*}
Indeed, let $p$ be the interpolating polynomial for $g$ constructed by the points of interpolation $\lambda_1$, \dots, $\lambda_{\m}$.
We take an arbitrary polynomial $q$ of degree at most $m$. Clearly the interpolating polynomial for $g-q$ is the polynomial $p-q$. Therefore, by~\eqref{e:norm(f-p)} we have
\begin{align*}
\max_{\lambda\in[\hat\alpha,\hat\beta]}|g(\lambda)-p(\lambda)|
&=\max_{\lambda\in[\hat\alpha,\hat\beta]}\bigl|\bigl(g(\lambda)-q(\lambda)\bigr)-\bigl(p(\lambda)-q(\lambda)\bigr)\bigr|\\
&\le(1+\Lambda_{\{\lambda_i\}})\max_{\lambda\in[\hat\alpha,\hat\beta]}|g(\lambda)-q(\lambda)|.
\end{align*}
Taking the minimum in this inequality over all $q$, we obtain the estimate being proved.

For the Lebesgue constant $\Lambda_{Ch}$ corresponding to $\m$ Chebyshev nodes, the following representation is known~\cite{Brutman97,Dzjadyk-Ivanov83}:
\begin{align*}
\Lambda_{Ch}&=\frac1\m\sum_{i=1}^\m\cot\frac{(2i-1)\pi}{4\m}.
\end{align*}
For 16 Chebyshev nodes, one has $\Lambda_{Ch}\approx2.69$. Thus, the Chebyshev nodes give almost the best uniform approximation.

We assume that the entries of the matrices $A$ and $T$ are known with the standard accuracy (IEEE double-precision binary floating-point format), which is about 16 significant decimal digits. Therefore, there is no need to know other parameters with greater accuracy if calculations involve $T$. In particular, knowing the interpolating polynomial (i.~e., the coefficients of the interpolating polynomial) with higher accuracy makes no sense. Furthermore, the approximation of the function $g$ by an interpolating polynomial $p$ with higher accuracy is also senseless. Moreover, the powers $T^k$ are calculated with a fewer number of significant digits than $T$ itself. Therefore, the leading coefficients of the interpolating polynomial may be calculated with lower accuracy.

Numerical experiments on the approximation of the function $\lambda\mapsto e^\lambda$ on the segment $[-\rho/2,\rho/2]$ by the interpolating polynomial $p$ of degree $m$ generated by the Chebyshev nodes are given in Table~\ref{tab:Ch}. One can see from Table~\ref{tab:Ch} that for $\rho=2$, the polynomial of the 14th degree provides about 16 significant decimal digits (for $\rho=1$, the polynomial of the 12th degree provides about 16 significant decimal digits). At the same time, the Paterson--Stockmeyer algorithm we use gives the best economy for the degree of a polynomial of the form $l^2-1$; the closest degree of the form $l^2-1$ to both the 14th and the 12th is the 15th. Therefore, we prefer $\rho=2$. We recall the definition of the Newton form of the interpolating polynomial below in formula~\eqref{e:poly Newton}.

\begin{table}
 \centering
 \caption{The value $\max_{\lambda\in[\hat\alpha,\hat\beta]}|e^\lambda-p_m(\lambda)|$ of the approximation of the function $\lambda\mapsto e^\lambda$ on $[\hat\alpha,\hat\beta]$ by the interpolating polynomial $p_m$ of degree $m$ with the scaled Chebyshev nodes taken as the points of interpolation}\label{tab:Ch}
\begin{tabular}{|c||c|c|c|c|c|c|}
  \hline
degree $m$ of $p_m$ & 7 & 13 & 14 & 15 \\
  \hline
  $\max_{\lambda\in[-1,1]}|e^\lambda-p_m(\lambda)|$ & $2.07\cdot10^{-7}$ & $1.46\cdot10^{-15}$ & $4.82\cdot10^{-17}$ & $1.51\cdot10^{-18}$ \\
  \hline
  \hline
degree $m$ of $p_m$ & 7 & 11 & 12 & 15 \\
  \hline
  $\max_{\lambda\in[-\frac12,\frac12]}|e^\lambda-p_m(\lambda)|$ & $8.02\cdot10^{-10}$ & $2.60\cdot10^{-16}$ & $4.98\cdot10^{-18}$ & $2.30\cdot10^{-23}$ \\
  \hline
\end{tabular}
\end{table}

Both the Lagrange and Newton forms of the interpolating polynomial contain differences $\lambda_i-\lambda_j$, where $\lambda_i$ are the points of interpolation. Therefore the approximate calculation of the interpolating polynomial accompanies by a loss of accuracy (for the Chebyshev nodes $\lambda_i$, the smallest of the numbers $|\lambda_i-\lambda_j|$ is about $0.098*(\beta_k-\alpha_k)$). We construct the interpolating polynomial in the Newton form because (according to numerical experiments) the Newton form gives 2 additional significant decimal digits in comparison with the Lagrange one. Since the calculation of the interpolating polynomial is a scalar (not matrix) operation, the usage of the enlarged number of decimal digits does not lead to an essential loss of time.

To achieve the accuracy in interpolating polynomials in 16 significant decimal digits, we propose calculating them with an enlarged number of decimal digits. Numerical experiments presented in Table~\ref{tab:Ch1} (with approximation of the exponential function by the interpolating polynomial with 16 Chebyshev nodes) shows that for $\beta_k-\alpha_k=2$, in order to provide a relative error of about $10^{-16}$ in all coefficients of the interpolating polynomial, it is necessary to use 33 significant decimal digits in the values of points of interpolation $\lambda_i$ (if $\beta_k-\alpha_k=1$, to provide a relative error of about $10^{-16}$ in coefficients, one must use 37 significant decimal digits\footnote{The `worse' accuracy for the smaller segment is explained by the fact that on a smaller segment the coefficients of the interpolating polynomial are closer to zero.}).

\begin{table}
 \centering
 \caption{The effect of the number of significant decimal digits (in the values of points of interpolation) on the error in the coefficients of the interpolating polynomial (of degree $15$ with the Chebyshev nodes taken as the points of interpolation) of the function $\lambda\mapsto e^\lambda$}\label{tab:Ch1}
\begin{tabular}{|c||c|c|c|c|c|}
  \hline
\multicolumn{4}{|c|}{on $[-1,1]$}\\
  \hline
significant decimal digits & 23 & 30 & 33 \\
  \hline
absolute error in the coefficients
 & $16.01$ & $23.01$ & $26.01$ \\
  \hline
relative error in the coefficients
 & $6.46$ & $13.46$ & $16.46$ \\
  \hline
 \hline
\multicolumn{4}{|c|}{on $[-1/2,1/2]$}\\
  \hline
significant decimal digits & 26 & 36 & 37 \\
  \hline
absolute error in the coefficients
 & $16.2$ & $26.2$ & $27.2$ \\
  \hline
relative error in the coefficients
 & $5.31$ & $15.3$ & $16.3$ \\
  \hline
\end{tabular}
\end{table}

Nevertheless, the largest coefficient (of the interpolating polynomial for the exponential function) is about 1 and the smallest one is about $7.76*10^{-13}$ for 16 Chebyshev nodes and $\rho=2$ (is about $7.67*10^{-13}$ for $\rho=1$). Table~\ref{tab:Ch1} shows that if one wants the absolute error in coefficients less than $10^{-16}$, it is enough to use 23 decimal digits in the case $\rho=2$ (and 26 decimal digits in the case $\rho=1$). In our numerical experiments (\S~\ref{s:numerical example}) we use 30 decimal digits in $\lambda_k$; such accuracy is needed for the calculation of the second block diagonal, see Section~\ref{s:DD}.

\begin{remark}\label{r:smale sizes}
If the order $N_k$ of the block $T_{kk}$ is small (e.g., $N_k<\m$), one can take for the points of interpolation the points of the spectrum of $T_{kk}$ and calculate $f(T_{kk})$ precisely, see~\cite[p.~5]{Higham08} for details.

We also recall that if $\beta_k-\alpha_k$ is small, one can use the Taylor polynomial instead of the interpolating one.

For an estimate of $\lVert p_{k}(A)-g(A)\rVert$, see~\cite{Kurbatov-Kurbatova-EMJ20}.
\end{remark}

\subsection{The usage of the Pad\'e approximation}\label{ss:Pade approximation}
Functions of matrices with a small spectrum are also often calculated by means of the Pad\'e approximation. In particular, if $f$ is the exponential function, the scaling and squaring method~\cite{Higham08} can be used. We do not discuss this way in this paper. We only note paper~\cite{Kressner-Luce-Statti17}, where the scaling and squaring method is applied to the successive calculation of the columns of a block triangular matrix.

\section{The Paterson--Stockmeyer algorithm}\label{s:Paterson-Stockmeyer}
Numerical experiments show that the accuracy provided by substituting a matrix directly into the Newton interpolating polynomial is insignificantly higher than the accuracy of substituting it into the same polynomial after expanding in powers of $\lambda$. Therefore, we substitute matrices into interpolating polynomials represented in the form
\begin{equation*}
p(\lambda)=c_0+c_1\lambda+c_2\lambda^2+c_3\lambda^3+c_4\lambda^4+c_5\lambda^5+\ldots+c_{m}\lambda^{m}
\end{equation*}
because when a polynomial is expanded in powers of $\lambda$, the Paterson--Stockmeyer
algorithm~\cite{Paterson-Stockmeyer73},~\cite[9.2.4]{Golub-Van_Loan13:eng} can be used. This algorithm allows one to reduce the number of matrix multiplications and thus speed up calculations. We describe here only one specific case of this algorithm that we use in our numerical experiments in Section~\ref{s:numerical example}. See also modifications of the Paterson--Stockmeyer algorithm in~\cite{Bader-Blanes-Casas19,Sastre18}.

Let $p$ be a polynomial of degree $15$ expanded in powers of $\lambda$:
\begin{equation*}
p(\lambda)=c_0+c_1\lambda+c_2\lambda^2+c_3\lambda^3+c_4\lambda^4+c_5\lambda^5+\ldots+c_{15}\lambda^{15}.
\end{equation*}
It is required to substitute into $p$ a given square matrix $A$, i.~e. to calculate
\begin{equation*}
p(A)=c_0\mathbf1+c_1A+c_2A^2+c_3A^3+c_4A^4+c_5A^5+\ldots+c_{15}A^{15}.
\end{equation*}

In this paper, by the Paterson--Stockmeyer algorithm we mean the following procedure.
First, we calculate the powers
\begin{equation*}
A^2,\; A^3,\; A^4.
\end{equation*}
This stage requires 3 matrix multiplications.

Next, we calculate the new coefficients (no matrix multiplications are needed)
\begin{align*}
c_0^{(1)}&=c_0\mathbf1+c_1A+c_2A^2+c_3A^3,\\
c_1^{(1)}&=c_4\mathbf1+c_5A+c_6A^2+c_7A^3,\\
c_2^{(1)}&=c_8\mathbf1+c_9A+c_{10}A^2+c_{11}A^3,\\
c_3^{(1)}&=c_{12}\mathbf1+c_{13}A+c_{14}A^2+c_{15}A^3.
\end{align*}
Now $p(A)$ can be represented in powers of $A^4$:
\begin{equation*}
p(A)=c_0^{(1)}+c_1^{(1)}A^4+c_2^{(1)}(A^4)^2+c_3^{(1)}(A^{4})^3.
\end{equation*}
For the final calculation of $p(A)$, we apply the Horner method:
\begin{equation*}
p(A)=c_0^{(1)}+\Bigl(c_1^{(1)}+\bigl(c_2^{(1)}+c_3^{(1)}A^4\bigr)A^4\Bigr)A^4.
\end{equation*}
The last calculation requires 3 additional matrix multiplications. In total, we need 6 matrix multiplications for calculating a polynomial of degree 15.

\section{Calculation of $F_{k,k+1}$}\label{s:DD}

Let $U\subset\mathbb{C}$ be an open set, and let $f:\,U\to\mathbb C$ be an analytic function. The \emph{divided difference} of $f$ is~\cite{Gelfond:eng,Jordan} the function $f^{[1]}:\,U\times
U\to\mathbb C$ defined by the formula
\begin{equation*}
f^{[1]}(\lambda,\mu)= \begin{cases}
\frac{f(\lambda)-f(\mu)}{\lambda-\mu}, & \text{if $\lambda\neq\mu$},\\
f'(\lambda), & \text{if $\lambda=\mu$}.
 \end{cases}
\end{equation*}

The simplest examples of divided differences are
\begin{align*}
v_1^{[1]}(\lambda,\mu)&=1&&\text{for }v_1(\lambda)=\lambda,\\
v_n^{[1]}(\lambda,\mu)&=\lambda^{n-1}+\lambda^{n-2}\mu+\dots+\mu^{n-1}&&\text{for }v_n(\lambda)=\lambda^n,\\
r_1^{[1]}(\lambda,\mu)&=\dfrac1{(\lambda_0-\lambda)(\lambda_0-\mu)}&&\text{for }r_1(\lambda)=\frac1{\lambda_0-\lambda}.
\end{align*}

\begin{proposition}[{\rm see, e.~g.,~\cite[Proposition 44]{Kurbatov-Kurbatova-Oreshina}}]\label{e:f[1] is analytic}
The function $f^{[1]}$ is analytic in $U\times U$.
\end{proposition}

The Taylor series for the divided difference of an analytic function $f$ at a point $(\lambda_0,\lambda_0)$ has the form
\begin{align*}
f^{[1]}(\lambda,\mu)&=\sum_{n=0}^\infty\frac{f^{(n+1)}(\lambda_0)}{(n+1)!}v_{n+1}^{[1]}(\lambda-\lambda_0,\mu-\lambda_0)\notag\\
&=\sum_{n=0}^\infty\frac{f^{(n+1)}(\lambda_0)}{(n+1)!}\sum_{i=0}^n(\lambda-\lambda_0)^{n-i}(\mu-\lambda_0)^i,
\end{align*}
where $v_n(\lambda)=\lambda^n$.
In particular, for $\exp_t(\lambda)=e^{\lambda t}$ at the point $(0,0)$ the expansion is
\begin{align*}
\exp_t^{[1]}(\lambda,\mu)&=\sum_{n=0}^\infty\frac{t^n}{(n+1)!}v_{n+1}^{[1]}(\lambda,\mu)
=\sum_{n=0}^\infty\frac{t^n}{(n+1)!}\sum_{i=0}^n\lambda^{n-i}\mu^i.
\end{align*}

Let $A$, $B$, and $H$ be matrices of the sizes $\nu\times\nu$, $\mu\times\mu$, and $\nu\times\mu$ respectively. We denote by $\sigma(A)$ and $\sigma(B)$ the spectra of $A$ and $B$. Let $f$ be an analytic function defined on an open set $U\subset\mathbb{C}$ such that $\sigma(A)\cup\sigma(B)\subset U$. We call the \emph{divided difference of $f$ applied to $A$ and $B$ at $H$} the matrix
\begin{equation}\label{e:f[1](A,B)}
f^{[1]}(A,B)\diamond H=\frac1{(2\pi
i)^2}\int_{\Gamma_2}\int_{\Gamma_1}f^{[1]}(\lambda,\mu)(\lambda\mathbf1-A)^{-1}H(\mu\mathbf1-B)^{-1}
\,d\lambda\,d\mu,
\end{equation}
where $\Gamma_1$ surrounds $\sigma(A)$ and $\Gamma_2$ surrounds $\sigma(B)$ in the counterclockwise direction, and $f$ is analytic inside both $\Gamma_1$ and $\Gamma_2$. For more about general functions of two matrices, see~\cite{Gil-MN08,Gil-EJLA-11,Kressner14,Kressner19,Kurbatov-Kurbatova-Oreshina}.

We mention the following statement (which we do not use explicitly).

\begin{proposition}[{\rm\cite[Theorem 45]{Kurbatov-Kurbatova-Oreshina}, \cite[Theorem 17]{Kurbatov-Kurbatova-OM19}}]\label{p:t:Theta}
Under the above assumptions
\begin{equation*}
f^{[1]}(A,B)\diamond H=\frac1{2\pi
i}\int_{\Gamma}f(\lambda)(\lambda\mathbf1-A)^{-1}H(\mu\mathbf1-B)^{-1}
\,d\lambda,
\end{equation*}
where $\Gamma$ surrounds $\sigma(A)\cup\sigma(B)$ in the counterclockwise direction.
\end{proposition}

\begin{proposition}[{\rm\cite{Carbonell-Jimenez-Pedroso08},~\cite{Davis_Ch73},~\cite{Kurbatov-Kurbatova-OM19}}]\label{p:Carbonell}
The blocks $F_{k,k+1}$ of the second diagonal of block matrix~\eqref{e:block form of F} can be represented in the form
\begin{equation}\label{e:F_{i,i+1}}
F_{k,k+1}=f^{[1]}(T_{kk},T_{k+1,k+1})\diamond T_{k,k+1}.
\end{equation}
\end{proposition}

In this section, we describe an algorithm for approximate calculation of~\eqref{e:F_{i,i+1}}.

We fix $k=1,2,\dots,n-1$.
We denote by $[\alpha_1,\beta_1]$ and $[\alpha_2,\beta_2]$ segments that contain the spectra $\sigma(T_{kk})$ and $\sigma(T_{k+1,k+1})$ respectively. We introduce the notation $\gamma_1=\frac{\alpha_1+\beta_1}2$ and $\gamma_2=\frac{\alpha_2+\beta_2}2$. We denote the Chebyshev nodes referred to the segments
\begin{equation*}
[\hat\alpha_1,\hat\beta_1]=[\alpha_1-\gamma_1,\beta_1-\gamma_1]\text{ and } [\hat\alpha_2,\hat\beta_2]=[\alpha_2-\gamma_2,\beta_2-\gamma_2]
\end{equation*}
by $\lambda_1$, $\lambda_2$, \dots, $\lambda_{\m_1}$ and
$\mu_1$, $\mu_2$, \dots, $\mu_{\m_2}$ respectively.

We consider the shifted function
\begin{equation*}
h(\lambda,\mu)=f^{[1]}(\lambda+\gamma_1,\mu+\gamma_{2}).
\end{equation*}
We set 
\begin{equation*}
A=T_{kk}-\gamma_1\mathbf1,\qquad B=T_{k+1,k+1}-\gamma_{2}\mathbf1,\qquad H=T_{k,k+1}.
\end{equation*}
By~\eqref{e:F_{i,i+1}} and the definition,
\begin{align*}
F_{k,k+1}&=f^{[1]}(T_{kk},T_{k+1,k+1})\diamond T_{k,k+1}\\
&=\frac1{(2\pi
i)^2}\int_{\Gamma_2}\int_{\Gamma_1}f^{[1]}(\lambda,\mu)(\lambda\mathbf1-T_{kk})^{-1}
T_{k,k+1}(\mu\mathbf1-T_{k+1,k+1})^{-1}
\,d\lambda\,d\mu\\
&=\frac1{(2\pi
i)^2}\int_{\Gamma_2}\int_{\Gamma_1}f^{[1]}(\lambda,\mu)(\lambda\mathbf1-T_{kk})^{-1}
H(\mu\mathbf1-T_{k+1,k+1})^{-1}
\,d\lambda\,d\mu.
\end{align*}
We make the change of variable $\lambda+\gamma_1\to\lambda$ in the internal integral:
\begin{align*}
\int_{\Gamma_1}f^{[1]}(\lambda,\mu)(\lambda\mathbf1-T_{kk})^{-1}
\,d\lambda
&=\int_{\Delta_1}f^{[1]}(\lambda+\gamma_1,\mu)(\lambda\mathbf1+\gamma_i\mathbf1-T_{kk})^{-1}
\,d\lambda\\
&=\int_{\Delta_1}f^{[1]}(\lambda+\gamma_1,\mu)(\lambda\mathbf1-A)^{-1}
\,d\lambda,
\end{align*}
where $\Delta_1=\Gamma_1-\gamma_1$ surrounds $\sigma(A)$.
We make the similar change of variable $\mu+\gamma_2\to\mu$ in the external integral:
\begin{multline*}
\int_{\Gamma_2}\int_{\Delta_1}f^{[1]}(\lambda+\gamma_1,\mu)(\lambda\mathbf1-A)^{-1}
H(\mu\mathbf1-T_{k+1,k+1})^{-1}
\,d\lambda\,d\mu\\
=\int_{\Delta_2}\int_{\Delta_1}f^{[1]}(\lambda+\gamma_1,\mu+\gamma_2)(\lambda\mathbf1-A)^{-1}
H(\mu\mathbf1+\gamma_2\mathbf1-T_{k+1,k+1})^{-1}
\,d\lambda\,d\mu\\
=\int_{\Delta_2}\int_{\Delta_1}f^{[1]}(\lambda+\gamma_1,\mu+\gamma_2)(\lambda\mathbf1-A)^{-1}
H(\mu\mathbf1-B)^{-1}
\,d\lambda\,d\mu\\
=\int_{\Delta_2}\int_{\Delta_1}h(\lambda,\mu)(\lambda\mathbf1-A)^{-1}
H(\mu\mathbf1-B)^{-1}
\,d\lambda\,d\mu,
\end{multline*}
where $\Delta_2=\Gamma_2-\gamma_2$ surrounds $\sigma(B)$. Thus,
\begin{equation}\label{e:F_{i,i+1}=h(A,B)H}
F_{k,k+1}=h(A,B)\diamond H.
\end{equation}

We construct an interpolating polynomial
\begin{equation}\label{e:p(lambda,mu)=:0}
p(\lambda,\mu)=\sum_{i=0}^{\m_1-1}\sum_{j=0}^{\m_2-1}c_{ij}\lambda^i\mu^j
\end{equation}
that agrees with $h$ at the points $(\lambda_i,\mu_j)$, i.~e.,
\begin{equation}\label{e:p=h}
p(\lambda_i,\mu_j)=h(\lambda_i,\mu_j),\qquad i=1,2,\dots,\m_1;\,j=1,2,\dots,\m_2.
\end{equation}
But first we recall some notation and terminology connected with the Newton interpolating polynomial.

Let be given \emph{distinct}\footnote{In this paper, we need not interpolating polynomials with multiple points of interpolation.} points of interpolation $z_1$, $z_2$, \dots,
$z_{\m}$ (complex numbers) and the family $f_1$, $f_2$, \dots,
$f_{\m}$ of complex numbers. The \emph{divided differences}~\cite{Gelfond:eng,Jordan} of the family $f_1$, $f_2$, \dots,
$f_{\m}$ (with respect to the points $z_1$, $z_2$, \dots,
$z_{\m}$) are the family $f_i^{i+j}$ defined by the recurrence relations
\begin{align*}
f_i^{i}&=f_i,\\
f_i^{i+1}&=\frac{f_{i+1}^{i+1}-f_i^{i}}{z_{i+1}-z_i},\\
\dots\dots&\dots\dots\dots\dots\dots\\
f_i^{i+j}&=\frac{f_{i+1}^{i+j}-f_i^{i+j-1}}{z_{i+j}-z_i}.
\end{align*}
It is convenient to arrange the divided differences in the table
\begin{equation}\label{e:table}
\begin{matrix}
f_1^{1}&f_2^{2}&\dots&f_i^{i}&\dots&f_{\m-2}^{\m-2}&f_{\m-1}^{\m-1}&f_{\m}^{\m}\\
f_1^{2}&f_2^{3}&\dots&f_i^{i+1}&\dots&f_{\m-2}^{\m-1}&f_{\m-1}^{\m}\\
f_1^{3}&f_2^{4}&\dots&f_i^{i+2}&\dots&f_{\m-2}^{\m}\\
\hdotsfor[2]{6}\\
f_1^{\m-1}&f_2^{\m}\\
f_1^{\m}
\end{matrix}
\end{equation}
We note that the construction of the Newton interpolating polynomial~\eqref{e:poly Newton} requires only the first column of this table. But to calculate the first column, all the others are needed.

We recall~\cite{Gelfond:eng,Jordan} that for any points of interpolation $z_1$, $z_2$, \dots, $z_{\m}$ and values $f_1$, $f_2$, \dots,
$f_{\m}$, the interpolating polynomial\footnote{An \emph{interpolation polynomial} is the polynomial $p$ of degree $m=\m-1$ satisfying the conditions $p(z_1)=f_1$, $p(z_2)=f_2$, \dots, $p(z_{\m})=f_{\m}$.} can be represented in the Newton form (also called the Newton interpolating polynomial)
\begin{equation}\label{e:poly Newton}
\begin{split}
p(z)&=f_1^{1}+f_1^{2}(z-z_1)
f_1^{3}(z-z_1)(z-z_2)\\
&+f_1^{4}(z-z_1)(z-z_2)(z-z_3)+\ldots\\
&+f_1^{\m-1}(z-z_1)(z-z_2)\ldots(z-z_{\m-2})\\
&+f_1^{\m}(z-z_1)(z-z_2)\ldots(z-z_{\m-1}).
\end{split}
\end{equation}

Now we begin the construction of~\eqref{e:p(lambda,mu)=:0}.

For each $j=1,2,\dots,\m_2$, we consider the function
\begin{equation*}
g_j(\lambda)=h(\lambda,\mu_j).
\end{equation*}
Then we construct the interpolating polynomial $b_j$ of degree $m=\m_1-1$ that agrees with the function $g_j$ at the points $\lambda_1$, $\lambda_2$, \dots, $\lambda_{\m_1}$, i.~e.,
\begin{equation}\label{e:q_i(lambda_j)}
b_j(\lambda_i)=g_j(\lambda_i)=h(\lambda_i,\mu_j),\qquad i=1,2,\dots,\m_1.
\end{equation}
We represent the interpolating polynomials $b_j$ in the Newton form
\begin{multline}\label{e:qj in the Newton form}
b_j(\lambda)=d_{0j}+d_{1j}(\lambda-\lambda_1)+d_{2j}(\lambda-\lambda_1)(\lambda-\lambda_2)\\
+\ldots+d_{\m_1-1,\,j}(\lambda-\lambda_1)\dots(\lambda-\lambda_{\m_1-1}),
\end{multline}
where $d_{0j}$, $d_{1j}$, \dots, $d_{\m_1-1,\,j}$ are the divided differences of the family $h(\lambda_1,\mu_j)$, $h(\lambda_2,\mu_j)$, \dots, $h(\lambda_{\m_1-1},\mu_j)$.

Then for each $i=0,1,\dots,\m_1-1$, we consider the family
\begin{equation*}
\begin{matrix}
(d_i)_1^{1}=d_{i,1},&(d_i)_2^{2}=d_{i,2},&\dots,&
&(d_i)_{\m_2-1}^{\m_2-1}=d_{i,\m_2-1},&(d_i)_{\m_2}^{\m_2}=d_{i,\m_2}
\end{matrix}
\end{equation*}
and calculate the divided differences with respect to the points of interpolation $\mu_1$, $\mu_2$, \dots, $\mu_{\m_2}$:
\begin{equation*}
\begin{matrix}
(d_i)_1^{1}&(d_i)_2^{2}&\dots&(d_i)_k^{k}&\dots&(d_i)_{\m_2-2}^{\m_2-2}&(d_i)_{\m_2-1}^{\m_2-1}&(d_i)_{\m_2}^{\m_2}\\
(d_i)_1^{2}&(d_i)_2^{3}&\dots&(d_i)_k^{k+1}&\dots&(d_i)_{\m_2-2}^{\m_2-1}&(d_i)_{\m_2-1}^{\m_2}\\
(d_i)_1^{3}&(d_i)_2^{4}&\dots&(d_i)_k^{k+2}&\dots&(d_i)_{\m_2-2}^{\m_2}\\
\hdotsfor[2]{6}\\
(d_i)_1^{\m_2-1}&(d_i)_2^{\m_2}\\
(d_i)_1^{\m_2}
\end{matrix}
\end{equation*}
After that we construct the Newton interpolating polynomial
\begin{equation}\label{e:ri}
r_i(\mu)=(d_i)_1^{1}+(d_i)_1^{2}(\mu-\mu_1)+\ldots+(d_i)_1^{\m_2}(\mu-\mu_1)\dots(\mu-\mu_{\m_2-1}).
\end{equation}
Since $r_i$ is the interpolating polynomial, one has
\begin{equation}\label{e:ri(muj)=d{ij}}
r_i(\mu_j)=d_{ij}.
\end{equation}

Finally, we consider the polynomial in two variables
\begin{multline}\label{e:p(lambda,mu):Newton}
p(\lambda,\mu)=r_0(\mu)+r_1(\mu)\;(\lambda-\lambda_1)+r_2(\mu)\;(\lambda-\lambda_1)(\lambda-\lambda_2)\\
+\ldots+r_{\m_1-1}(\mu)\;(\lambda-\lambda_1)\dots(\lambda-\lambda_{\m_1-1}).
\end{multline}

From~\eqref{e:ri(muj)=d{ij}} one has
\begin{multline*}
p(\lambda,\mu_j)=d_{0j}+d_{1j}(\lambda-\lambda_1)+d_{2j}(\lambda-\lambda_1)(\lambda-\lambda_2)\\
+\ldots+d_{\m_1-1,\,j}(\lambda-\lambda_1)\dots(\lambda-\lambda_{\m_1-1})=b_j(\lambda).
\end{multline*}
In particular, from~\eqref{e:q_i(lambda_j)} it follows that conditions~\eqref{e:p=h} are fulfilled:
\begin{equation*}
p(\lambda_i,\mu_j)=b_j(\lambda_i)=h(\lambda_i,\mu_j),\qquad i=1,2,\dots,\m_1;\,j=1,2,\dots,\m_2.
\end{equation*}
Thus, the interpolating polynomial is constructed.
Expanding the parentheses in~\eqref{e:p(lambda,mu):Newton} and~\eqref{e:ri} we arrive at representation~\eqref{e:p(lambda,mu)=:0}. For more about interpolation of functions of several variables, see~\cite[ch.~5]{Phillips03}.

Finally, according to~\eqref{e:F_{i,i+1}=h(A,B)H}, we arrive at the approximation
\begin{equation}\label{e:F_{i,i+1} approx}
F_{k,k+1}=h(A,B)\diamond H\approx p(A,B)\diamond H=
\sum_{i=0}^{\m_1-1}\sum_{j=0}^{\m_2-1}c_{ij}A^iHB^j.
\end{equation}

\begin{remark}\label{r:p(A,B):Newton}
Clearly, the same approximation is obtained if one sets (according to~\eqref{e:p(lambda,mu):Newton})
\begin{multline*}
h(A,B)\diamond H\approx Hr_0(B)+(A-\lambda_1\mathbf1)Hr_1(B)+(A-\lambda_1\mathbf1)(A-\lambda_2\mathbf1)Hr_2(B)\\
+\ldots+(A-\lambda_1\mathbf1)\dots(A-\lambda_{\m_1-1}\mathbf1)Hr_{\m_1-1}(B).
\end{multline*}
We do not use this representation because we want to apply the Paterson--Stockmeyer algorithm.
\end{remark}

\begin{remark}\label{r:p in unique}
One can interpret equalities~\eqref{e:p=h} as conditions on the coefficients $c_{ij}$ of polynomial~\eqref{e:p(lambda,mu)=:0}. Thus he obtains a system of $\m_1\cdot \m_2$ linear equations with $\m_1\cdot \m_2$ unknowns. Above, we have proved that this system has a solution for any free terms $h(\lambda_i,\mu_j)$. Consequently, this solution (which uniquely determines the interpolating polynomial) is unique.
\end{remark}

The usage of the formula
\begin{equation*}
f^{[1]}(\lambda_i,\mu_j)
=\frac{f(\lambda_i)-f(\mu_j)}{\lambda_i-\mu_j}
\end{equation*}
can cause a loss of accuracy if $\lambda_i$ and $\mu_j$ are close\footnote{If one uses 16 Chebyshev nodes as points of interpolation and $\rho=2$, then the minimal value of $|\lambda_i-\mu_j|$ is about $0.0096$ (for $\rho=1$, the minimal value of $|\lambda_i-\mu_j|$ is about $0.0048$).}. Numerical experiments presented in Table~\ref{tab:Ch2} (with 16 Chebyshev nodes and the exponential function) show that 30 significant decimal digits provide the absolute error in $h(\lambda,\mu)$ about $10^{-16}$ for $(\lambda,\mu)\in\sigma(A)\times\sigma(B)$ in the case $\rho=2$ (but 36 decimal digits provide the absolute error in $h(\lambda,\mu)$ about $10^{-16}$ in the case $\rho=1$),
and 59 significant decimal digits provide the relative error in $h(\lambda_i,\mu_j)$ about $10^{-16}$ for the case $\rho=2$
(but 68 decimal digits provide the relative error in $h(\lambda_i,\mu_j)$ about $10^{-16}$ for the case $\rho=1$). We prefer to ensure the absolute (as opposed to reative) error for the reasons described in Section~\ref{ss:interpolation polynomial}.

\begin{table}
 \centering
 \caption{The effect of the number of significant decimal digits (in the values of points of interpolation) on the error in the coefficients of the interpolating polynomial (in two variables of degree $(15,15)$ with the Chebyshev nodes taken as the points of interpolation) of the function $\exp^{[1]}(\lambda,\mu)=\frac{e^\lambda-e^\mu}{\lambda-\mu}$}\label{tab:Ch2}
\begin{tabular}{|c||c|c|c|c|c|}
  \hline
\multicolumn{4}{|c|}{on $[-3,-1]\times[-1,1]$}\\
  \hline
significant decimal digits & 30 & 44 & 59 \\
  \hline
absolute error in the coefficients
 & $16.01$ & $30.10$ & $47.8$ \\
  \hline
relative error in the coefficients
 & $0$ & $1.08$ & $16.16$ \\
  \hline
 \hline
\multicolumn{4}{|c|}{on $[-3/2,-1/2]\times[-1/2,1/2]$}\\
  \hline
significant decimal digits & 36 & 53 & 68 \\
  \hline
absolute error in the coefficients
 & $15.8$ & $32.8$ & $45.10$ \\
  \hline
relative error in the coefficients
 & $0$ & $1.16$ & $16.01$ \\
  \hline
\end{tabular}
\end{table}

We are now in a position to describe the algorithm of approximate calculation of $h(A,B)\diamond H$ in full. For definiteness, let $\rho=2$.
We calculate the points of interpolation $\lambda_1$, $\lambda_2$, \dots, $\lambda_{\m_1}$ and
$\mu_1$, $\mu_2$, \dots, $\mu_{\m_2}$ with 30 significant decimal digits (actually, they should be already calculated when calculating the diagonal blocks $f(T_{kk})$ and $f(T_{k+1,k+1})$). Then we calculate the numbers
\begin{equation*}
p(\lambda_i,\mu_j)=b_j(\lambda_i)=h(\lambda_i,\mu_j),\qquad i=1,2,\dots,\m_1;\,j=1,2,\dots,\m_2.
\end{equation*}
After that, for each $j=1,2,\dots,\m_2$, we find the Newton interpolating polynomials $b_j$ that agrees with the function $\lambda\mapsto h(\lambda,\mu_j)$ at the points $\lambda_1$, $\lambda_2$, \dots, $\lambda_{\m_1}$, see formula~\eqref{e:qj in the Newton form}.

Then, for each $i=1,2,\dots,\m_1$, we calculate polynomials~\eqref{e:ri}, substitute them in~\eqref{e:p(lambda,mu):Newton}, and expand the parenthesis. As a result, we obtain the approximating polynomial $p$ in the form~\eqref{e:p(lambda,mu)=:0}. It remains to substitute $A$, $B$, and $H$ into it and obtain approximation~\eqref{e:F_{i,i+1} approx}:
\begin{equation*}
h(A,B)\diamond H\approx\sum_{j=0}^{\m_2-1}\sum_{i=0}^{\m_1-1}c_{ij}A^iHB^j.
\end{equation*}
To this end, for each $j=1,2,\dots,\m_2$, we calculate the sums
\begin{equation*}
\sum_{i=0}^{\m_1-1}c_{ij}A^i.
\end{equation*}
To accelerate the process, we apply the Paterson--Stockmeyer algorithm. Then we multiply the results by $H$ and obtain the matrices
\begin{equation*}
C_j=\sum_{i=0}^{\m_1-1}c_{ij}A^iH,\qquad j=1,2,\dots,\m_2.
\end{equation*}
It remains to calculate the sum
\begin{equation*}
\sum_{j=0}^{\m_2-1}C_jB^j.
\end{equation*}
Now the direct application of the Paterson--Stockmeyer algorithm requires 18 matrix multiplications for $\m_2=16$ because the coefficients $C_j$ are matrices. Therefore we use the Horner method~\cite[\S~4.2]{Higham08}, which requires $\m_2-1$ matrix multiplications.

 \begin{remark}\label{r:Kenney-Laub98:tanh}
In~\cite{McCurdy-Ng-Parlett84} and~\cite[formula (10.17)]{Higham08}, it is proposed the representations
\begin{equation*}
\exp_t^{[1]}(\lambda,\mu)=e^{\frac{(\lambda+\mu)t}2}\frac{\sinh\bigl(\frac{\lambda-\mu}2t\bigr)}{\frac{\lambda-\mu}2t},\qquad
\exp_t^{[1]}(\lambda,\mu)=e^{\mu t}\frac{e^{(\lambda-\mu)t}-1}{\lambda-\mu},\qquad\lambda\neq\mu,
\end{equation*}
for $\exp_t(\lambda)=e^{\lambda t}$.
In~\cite{Kenney-Laub98} another representation of a similar kind is proposed:
\begin{equation*}
\exp_t^{[1]}(\lambda,\mu)=(e^{\lambda t}+e^{\mu t})\frac{\tanh
\bigl(\frac{\lambda-\mu}{2}t\bigr)}{\frac{\lambda-\mu}{2}t},\qquad\lambda\neq\mu.
\end{equation*}
These representations can be used for the calculation of $\exp_t^{[1]}(A,B)\diamond H$. For example, we describe the idea for the case of the first formula. The function
\begin{equation*}
\sinch(z)=\frac{\sinh z}{z}
\end{equation*}
is entire (i.~e. analytic in the whole $\mathbb C$). Therefore its Taylor series
\begin{equation*}
\sinch(z)=\sum_{k=0}^\infty\frac{z^{2k}}{(2k+1)!}=1+\frac{z^2}{3!}+\frac{z^4}{5!}+\dots.
\end{equation*}
converges everywhere. Hence $\exp_t^{[1]}(A,B)\diamond H$ can be calculated by the application of the function $z\mapsto\sinch(z)$ to the transformator $H\mapsto(AH-HB)t/2$:
\begin{align*}
H\mapsto H+\frac{t^2}{2^2\cdot3!}&\Bigl(A\bigl(A(AH-HB)-(AH-HB)B\bigr)\\
&-\bigl(A(AH-HB)-(AH-HB)B\bigr)B\Bigr)+\dots
\end{align*}
and then the application to the result the transformator $G\mapsto e^{At/2}Ge^{Bt/2}$:
\begin{multline*}
\exp_t^{[1]}(A,B)\diamond H =e^{At/2}\Bigl(H+\frac{t^2}{2^2\cdot3!}\Bigl(A\bigl(A(AH-HB)-(AH-HB)B\bigr)\\
-\bigl(A(AH-HB)-(AH-HB)B\bigr)B\Bigr)+\dots\Bigr)e^{Bt/2}.
\end{multline*}
See~\cite[\S~10.2]{Higham08} or~\cite[Remark 1]{Kurbatov-Kurbatova-OM19} for details.
\end{remark}

We note the papers~\cite{Al-Mohy-Higham09,Al-Mohy-Higham10,Al-Mohy-Higham-Relton13}, in which the problem of calculating $f^{[1]}(A,A)\diamond H$ for $A=B$ is discussed in great detail.

\section{The whole algorithm}\label{s:algorithm}
In this section, we recall the sequence of our manipulations.

0. Given a square matrix $R$ with a real spectrum and a function $f$ analytic in an open neighborhood of the spectrum of $R$.
We assume that the entries of $R$ are known with the standard accuracy (IEEE double-precision binary floating-point format, which admits about 16 significant decimal digits) and $f$ can be calculated with any desirable accuracy. The aim is the approximate calculation of $f(R)$.

1. We construct the Schur decomposition $R=QTQ^H$, where $Q$ is unitary and $T$ is upper triangular. We denote by $F$ the matrix $f(T)$.

2. We find a  half-open interval\footnote{We recall that one can use a closed interval $[a,b]$. In this case, the notation becomes less symmetric but nothing changes essentially.} $[a,b)$ that contains the spectrum $\sigma(T)$ of $T$ (we recall that $\sigma(T)$ coincides with the set of diagonal entries of $T$).

3. We choose a number $\rho>0$. In our numerical experiments in Section~\ref{s:numerical example} we set $\rho=2$. Then we divide $[a,b)$ into parts of the length $\rho$:
\begin{equation*}
[a,b)\subseteq[a,a+\rho)\cup[a+\rho,a+2\rho)
\cup\dots\cup[a+(n-1)\rho,a+n\rho).
\end{equation*}
Here $n$ denotes the number of the parts.

4. We split the spectrum of $T$ into $n$ clusters: the $k$-th cluster is contained in $[a+k\rho,a+(k+1)\rho)$. If some clusters are empty, we skip them and enumerate the rest anew. Then we reorder (see Section~\ref{s:reodering}) the diagonal entries of $T$ so that the entries of the first cluster precede those of the second, the entries of the second cluster precede those of the third, and so on. The order of entries inside a cluster does not matter. Thus, we come to the matrix $T'=UTU^T$, where $U$ is the product of unitary similarity transformations that reorder the neighboring diagonal entries, see Section~\ref{s:reodering}. Below we denote $T'$ by the original symbol $T$.

5. For each $k=1,2,\dots,n$, we calculate $\alpha_k$, the minimal eigenvalue of $T_{kk}$, and $\beta_k$, the maximal one. We transform $\alpha_k$ and $\beta_k$ into numbers having a large number of significant digits (e.~g., by appending additional digits after the already available reliable ones and rounding up or down), which will ensure 16 significant digits in the coefficients of $p$ in step 7. In our numerical experiments in Section~\ref{s:numerical example}, we use 30 significant digits in $\alpha_k$ and $\beta_k$.
We calculate $\gamma_k=\frac{\alpha_k+\beta_k}2$, the center of the segment $[\alpha_k,\beta_k]$, and its radius $\delta_k=\frac{\beta_k-\alpha_k}2$.

6. We calculate the Chebyshev nodes $z_1,z_2\dots,z_{16}$ on the segment $[-1,1]$ with a large number of significant decimal digits (30 in our numerical experiments in Section~\ref{s:numerical example}).
For each $k=1,2,\dots,n$, we calculate the Chebyshev nodes $\lambda_1$, $\lambda_2$, \dots, $\lambda_{16}$ on the segment $[\hat\alpha,\hat\beta]=[-\delta_k,\delta_k]$ by the scaling
$\lambda_k=\delta_kz_k$; thus, the numbers $\lambda_k$ have about 30 significant decimal digits. Then we construct the Newton interpolating polynomial $p$ that agrees with the function $g(\lambda)=f(\lambda+\gamma_k)$ at $\lambda_1$, $\lambda_2$, \dots, $\lambda_{16}$. We expand $p$ in powers of $\lambda$. Using the Paterson--Stockmeyer algorithm (see \S~\ref{s:Paterson-Stockmeyer}) we calculate the approximation $F_{kk}=g(A)\approx p(A)$, where $A=T_{kk}-\gamma_k\mathbf1$.

7. For each $k=1,2,\dots,n-1$, we consider the shifted function $h(\lambda,\mu)=f^{[1]}(\lambda+\gamma_k,\mu+\gamma_{k+1})$, where $f^{[1]}(\lambda,\mu)=
\frac{f(\lambda)-f(\mu)}{\lambda-\mu}$. We denote by $\lambda_1$, $\lambda_2$, \dots, $\lambda_{16}$ and $\mu_1$, $\mu_2$, \dots, $\mu_{16}$ the Chebyshev nodes on $[\hat\alpha_1,\hat\beta_1]=[\alpha_k-\gamma_k,\beta_k-\gamma_k]$ and $[\hat\alpha_{2},\hat\beta_{2}]=[\alpha_{k+1}-\gamma_{k+1},\beta_{k+1}-\gamma_{k+1}]$ respectively (calculated in step 6).
We construct the interpolating polynomial $p$ in two variables that agrees with the function $h$ at the points $(\lambda_i,\mu_j)$ as it is described in Section~\ref{s:DD}. (To ensure higher accuracy, we first construct the interpolating polynomial $p$ in the Newton form~\eqref{e:p(lambda,mu):Newton}, and then expand it in powers of $\lambda$ and $\mu$.) We substitute the matrices $A=T_{kk}-\gamma_k\mathbf1$, $B=T_{k+1,k+1}-\gamma_{k+1}\mathbf1$, and $H=T_{k,k+1}$ into $p$ (see formula~\eqref{e:F_{i,i+1} approx}) and obtain the approximation $F_{k,k+1}=h(A,B)\diamond H\approx p(A,B)\diamond H$.

8. We fill in the remaining overdiagonal entries of $F$ by formula~\eqref{e:Parlett:single} (diagonal by diagonal). Finally, we return to the matrix $R$ by means of $f(R)=QUFU^HQ^H$.

\section{Numerical experiments}\label{s:numerical example}
We performed numerical experiments using `Mathematica'~\cite{Wolfram}.

For the order of the matrix $T$ we take the number $N=256$. We set $n=5$ (the number of blocks). We take four integers $w_1$, $w_2$, $w_3$, and $w_4$ generated by a random variable uniformly distributed in $[0,N]$. We arrange $w_1$, $w_2$, $w_3$, and $w_4$ in the increasing order. We set $w_0=0$, $w_5=N$ and define the orders $N_k$ of the blocks by the formula\footnote{If $N_k<2$ for some $k$ (this case is less interesting), we create $w_1$, $w_2$, $w_3$, and $w_4$ all over again.} $N_k=w_k-w_{k-1}$, $k=1,2,\dots,n$. We create the matrix $D$ of the size $N\times N$ consisting of zeroes and represent $D$ in the block form. For each $k=1,2,\dots,n$, we replace the block $D_{kk}$ by a diagonal matrix with diagonal entries uniformly distributed in $[-2k,-2k+2)$ and having 128 significant digits. We create an upper triangular matrix $S$, consisting of units on the main diagonal and random numbers with 128 significant digits uniformly distributed in $[-1,1]$ above the main diagonal. We set
\begin{equation*}
T=SDS^{-1}.
\end{equation*}
Thus, we obtain a random triangular (block) matrix $T$ with diagonal entries arranged in the desired (i.~e., as after step 4 in Section~\ref{s:algorithm}) order.

We take for $f$ the function $f(\lambda)=e^\lambda$. For the precise matrix $\widetilde{F}=f(T)$ we take $Sf(D)S^{-1}$. After that, we round the matrices $T$ and $\widetilde{F}$ to 16 significant decimal digits.

We take $\rho=2$ and divide the half-open interval $[-10,0)$ (that contains the spectrum $\sigma(T)$ of $T$) into 5 parts:
\begin{equation*}
[-10,0)=[-10,-8)\cup[-8,-6)\cup[-6,-4)\cup[-4,-2)\cup[-2,0).
\end{equation*}
which generates the division of the spectrum $\sigma(T)$ into clusters. Clearly, the block structure of $T$ is in agreement with this division.
We calculate $F=f(T)$ according to the algorithm described in items 5--7 of Section~\ref{s:algorithm}.

We repeat this experiment 100 times. The result is as follows. The mean value\footnote{The mean value of numbers $x_1$, $x_2$, \dots, $x_{100}$ is defined to be
$\frac1{100}\sum_{i=1}^{100}x_i$.} of the condition number
\begin{equation*}
\varkappa(S)=\lVert S\rVert\cdot\lVert S^{-1}\rVert
\end{equation*}
is $2.11\cdot10^{16}$. The mean value of the condition number
\begin{equation*}
\varkappa(T)=\lVert T\rVert\cdot\lVert T^{-1}\rVert
\end{equation*}
is $5.76\cdot10^{27}$.
The mean value of the max-norm
\begin{equation*}
\lVert\widetilde{F}\rVert_\infty=\max_{i,j=1,\dots,N}|\tilde f_{ij}|
\end{equation*}
of $\widetilde{F}$ is $4.80\cdot10^9$.
The mean value of the operator norm $\lVert F-\widetilde{F}\rVert_{2\to2}$ of the error $F-\widetilde{F}$ is $0.0035$; here $\lVert\cdot\rVert_{2\to2}$ is the matrix norm induced by the Euclidean norm on $\mathbb C^N$.
The mean value of the Frobenius norm
\begin{equation*}
\lVert F-\widetilde{F}\rVert_2=\sqrt{\sum_{i=1}^N\sum_{j=1}^N\bigl| f_{ij}-\tilde f_{ij}\bigr|^2}
\end{equation*}
is $0.0045$.
The mean value of the incomplete Frobenius norm
\begin{equation*}
\lVert F-\widetilde{F}\rVert_{2,\;\text{incomplete}}=\sqrt{\sum\bigl| f_{ij}-\tilde f_{ij}\bigr|^2}
\end{equation*}
is $1.86\cdot10^{-10}$; here the sum is taken over only the entries of the two main block diagonals. The mean value of the max-norm
\begin{equation*}
\lVert F-\widetilde{F}\rVert_\infty=\max_{i,j=1,\dots,N}|f_{ij}-\tilde f_{ij}|
\end{equation*}
is $0.0020$; here the sum is taken over all $i,j=1,2,\dots,N$. The mean value of the incomplete max-norm
\begin{equation*}
\lVert F-\widetilde{F}\rVert_{\infty,\;\text{incomplete}}=\max|f_{ij}-\tilde f_{ij}|
\end{equation*}
is $1.014\cdot10^{-10}$; here the maximum is taken over the entries of only the two main block diagonals.
The mean value of the max-relative error
\begin{equation*}
r_\infty(F,\widetilde{F})=\max\frac{|f_{ij}-\tilde f_{ij}|}{|\tilde f_{ij}|},
\end{equation*}
is $8.54\cdot10^{-10}$; here the maximum is taken over all $i$ and $j$ such that $\tilde f_{ij}\neq0$.

We note that a small number of entries of the matrix $F$ (they are located in the upper right corner) are significantly larger in absolute value than the rest.
As a result, the absolute error is almost independent of the type of the norm (Frobenius or operator).

We repeat the above experiments with smaller condition numbers but a more complicated function $f$. Now we take the entries of $S$ from $[-1/2,1/2]$ instead of $[-1,1]$ and the function $f(\lambda)=\cos\lambda$. The result is as follows. The mean value of the condition number $\varkappa(S)$ is $5.14\cdot10^{5}$. The mean value of the condition number $\varkappa(T)$ is $9.94\cdot10^{7}$.
The mean value of the max-norm $\lVert\widetilde{F}\rVert_\infty$ is $4.3\cdot10^3$.
The mean value of the operator norm $\lVert F-\widetilde{F}\rVert_{2\to2}$ is $0.0017$.
The mean value of the Frobenius norm
$\lVert F-\widetilde{F}\rVert_2$ is $0.0018$.
The mean value of the incomplete Frobenius norm
$\lVert F-\widetilde{F}\rVert_{2,\;\text{incomplete}}$
is $1.5\cdot10^{-5}$. The mean value of the max-norm
$\lVert F-\widetilde{F}\rVert_\infty$
is $5.5\cdot10^{-4}$. The mean value of the incomplete max-norm
$\lVert F-\widetilde{F}\rVert_{\infty,\;\text{incomplete}}$ is $3.9\cdot10^{-6}$.
The mean value of the max-relative error is $r_\infty(F,\widetilde{F})=0.0074$.


\providecommand{\bysame}{\leavevmode\hbox to3em{\hrulefill}\thinspace}
\providecommand{\MR}{\relax\ifhmode\unskip\space\fi MR }
\providecommand{\MRhref}[2]{%
  \href{http://www.ams.org/mathscinet-getitem?mr=#1}{#2}
}
\providecommand{\href}[2]{#2}


\begin{thebibliography}{10}

\bibitem{Al-Mohy-Higham09}
A.~H. Al-Mohy and N.~J. Higham, \emph{Computing the {F}r\'{e}chet derivative of
  the matrix exponential, with an application to condition number estimation},
  SIAM J. Matrix Anal. Appl. \textbf{30} (2008/09), no.~4, 1639--1657.
  \MR{2486857}

\bibitem{Al-Mohy-Higham10}
\bysame, \emph{The complex step approximation to the {F}r\'{e}chet derivative
  of a matrix function}, Numer. Algorithms \textbf{53} (2010), no.~1, 113--148.
  \MR{2566131}

\bibitem{Al-Mohy-Higham-Relton13}
A.~H. Al-Mohy, N.~J. Higham, and S.~D. Relton, \emph{Computing the
  {F}r\'{e}chet derivative of the matrix logarithm and estimating the condition
  number}, SIAM J. Sci. Comput. \textbf{35} (2013), no.~4, C394--C410.
  \MR{3080997}

\bibitem{Babenko:eng:2nd}
K.~I. Babenko, \emph{{Osnovy chislennogo analiza} [{F}undamentals of numerical
  analysis]}, second ed., Regulyarnaya i {K}haoticheskaya {D}inamika,
  Moscow--Izhevsk, 2002, (in Russian).

\bibitem{Bader-Blanes-Casas19}
P.~Bader, S.~Blanes, and Casas F., \emph{Computing the matrix exponential with
  an optimized {T}aylor polynomial approximation}, Mathematics \textbf{7}
  (2019), no.~12, 1174.

\bibitem{Bai-Demmel98}
Zh. Bai and J.~W. Demmel, \emph{Using the matrix sign function to compute
  invariant subspaces}, SIAM J. Matrix Anal. Appl. \textbf{19} (1998), no.~1,
  205--225. \MR{1609964}

\bibitem{Bai-Demmel93}
\bysame, \emph{On swapping diagonal blocks in real {S}chur
  form}, Linear Algebra Appl. \textbf{186} (1993), 73--95. \MR{1217200}

\bibitem{Bernstein1914}
S.~Bernstein, \emph{Sur la meilleure approximation de {$|x|$} par des polynomes
  de degr\'{e}s donn\'{e}s}, Acta Math. \textbf{37} (1914), no.~1, 1--57, (in
  French). \MR{1555093}

\bibitem{Brutman97}
L.~Brutman, \emph{Lebesgue functions for polynomial interpolation---a survey},
  Ann. Numer. Math. \textbf{4} (1997), no.~1-4, 111--127. \MR{1422674}

\bibitem{Burrage-Hale-Kay12}
K.~Burrage, N.~Hale, and D.~Kay, \emph{An efficient implicit {FEM} scheme for
  fractional-in-space reaction-diffusion equations}, SIAM J. Sci. Comput.
  \textbf{34} (2012), no.~4, A2145--A2172. \MR{2970400}

\bibitem{Carbonell-Jimenez-Pedroso08}
F.~Carbonell, J.~C. J\'imenez, and L.~M. Pedroso, \emph{Computing multiple
  integrals involving matrix exponentials}, J. Comput. Appl. Math. \textbf{213}
  (2008), no.~1, 300--305. \MR{2382698}

\bibitem{Cardoso-Sadeghi18}
J.~R. Cardoso and A.~Sadeghi, \emph{Computation of matrix gamma function}, BIT
  \textbf{59} (2019), no.~2, 343--370. \MR{3974043}

\bibitem{Davies-Higham03}
Ph.~I. Davies and N.~J. Higham, \emph{A {S}chur-{P}arlett algorithm for
  computing matrix functions}, SIAM J. Matrix Anal. Appl. \textbf{25} (2003),
  no.~2, 464--485 (electronic). \MR{2047429}

\bibitem{Davis_Ch73}
Ch. Davis, \emph{Explicit functional calculus}, Linear Algebra and Appl.
  \textbf{6} (1973), 193--199. \MR{0327792}

\bibitem{Dzjadyk-Ivanov83}
V.~K. Dzjadyk and V.~V. Ivanov, \emph{On asymptotics and estimates for the
  uniform norms of the {L}agrange interpolation polynomials corresponding to
  the {C}hebyshev nodal points}, Anal. Math. \textbf{9} (1983), no.~2, 85--97.
  \MR{720078}

\bibitem{Gelfond:eng}
A.~O. Gel$'$fond, \emph{Calculus of finite differences}, International
  Monographs on Advanced Mathematics and Physics, Hindustan Publishing Corp.,
  Delhi, 1971, Translation of the third Russian edition. \MR{0342890}

\bibitem{Gil-MN08}
M.~I. Gil', \emph{Norm estimates for functions of two operators on tensor
  products of {H}ilbert spaces}, Math. Nachr. \textbf{281} (2008), no.~8,
  1129--1141. \MR{2427165}

\bibitem{Gil-EJLA-11}
\bysame, \emph{Norm estimates for functions of two non-commuting matrices},
  Electron. J. Linear Algebra \textbf{22} (2011), 504--512.

\bibitem{Godunov94:ODE:eng}
S.~K. Godunov, \emph{Ordinary differential equations with constant
  coefficient}, Translations of Mathematical Monographs, vol. 169, American
  Mathematical Society, Providence, RI, 1997, Translated from the 1994 Russian
  original. \MR{1465434}

\bibitem{Golub-Van_Loan13:eng}
G.~H. Golub and Ch.~F. Van~Loan, \emph{Matrix computations}, fourth ed., Johns
  Hopkins Studies in the Mathematical Sciences, Johns Hopkins University Press,
  Baltimore, MD, 2013. \MR{3024913}

\bibitem{Grimm-Hochbruck08}
V.~Grimm and M.~Hochbruck, \emph{Rational approximation to trigonometric
  operators}, BIT \textbf{48} (2008), no.~2, 215--229. \MR{2430617}

\bibitem{Hargreaves-Higham05}
G.~I. Hargreaves and N.~J. Higham, \emph{Efficient algorithms for the matrix
  cosine and sine}, Numer. Algorithms \textbf{40} (2005), no.~4, 383--400.
  \MR{2191973}

\bibitem{Higham08}
N.~J. Higham, \emph{Functions of matrices: theory and computation}, Society for
  Industrial and Applied Mathematics (SIAM), Philadelphia, PA, 2008.
  \MR{2396439}

\bibitem{Higham-Lin2011}
N.~J. Higham and L.~Lin, \emph{A {S}chur-{P}ad\'e algorithm for fractional
  powers of a matrix}, SIAM J. Matrix Anal. Appl. \textbf{32} (2011), no.~3,
  1056--1078.

\bibitem{Householder64}
A.~S. Householder, \emph{The theory of matrices in numerical analysis},
  Blaisdell Publishing Co. Ginn and Co., New York--Toronto--London, 1964.
  \MR{0175290}

\bibitem{Ikramov84:rus-eng}
Kh.~D. Ikramov, \emph{Numerical solution of matrix equations. {O}rthogonal
  methods}, Izdat. ``Nauka'', Moscow, 1984, (in Russian). \MR{779465}

\bibitem{Ikramov91a:eng}
Kh.~D. Ikramov, \emph{Matrix pencils --- theory, applications, numerical
  methods}, Mathematical analysis, {V}ol.\ 29, Itogi Nauki i Tekhniki, Akad.
  Nauk SSSR, Vsesoyuz. Inst. Nauchn. i Tekhn. Inform., Moscow, 1991, (in
  Russian); English translation in {\em J. Soviet Math.},
  \textbf{64}(2):783--853, 1993, pp.~3--106. \MR{1138216}

\bibitem{Jarlebring07}
E.~Jarlebring and T.~Damm, \emph{The {Lambert $W$} function and the spectrum of
  some multidimensional time-delay systems}, Automatica \textbf{43} (2007),
  no.~12, 2124--2128.

\bibitem{Jodar-Cortes}
L.~J\'odar and J.~C. Cort\'es, \emph{Some properties of gamma and beta matrix
  functions}, Appl. Math. Lett. \textbf{11} (1998), no.~1, 89--93. \MR{1490386}

\bibitem{Jordan}
Ch. Jordan, \emph{Calculus of finite differences}, third ed., Chelsea
  Publishing Co., New~York, 1965. \MR{0183987}

\bibitem{Kenney-Laub95}
C.~S. Kenney and A.~J. Laub, \emph{The matrix sign function}, IEEE Trans.
  Automat. Control \textbf{40} (1995), no.~8, 1330--1348. \MR{1343800}

\bibitem{Kenney-Laub98}
\bysame, \emph{A {S}chur--{F}r\'echet algorithm for computing the logarithm and
  exponential of a matrix}, SIAM J. Matrix Anal. Appl. \textbf{19} (1998),
  no.~3, 640--663 (electronic). \MR{1611163}

\bibitem{Kressner14}
D.~Kressner, \emph{Bivariate matrix functions}, Oper. Matrices \textbf{8}
  (2014), no.~2, 449--466. \MR{3224819}

\bibitem{Kressner19}
\bysame, \emph{A {K}rylov subspace method for the approximation of bivariate
  matrix functions}, Structured matrices in numerical linear algebra, Springer
  INdAM Ser., vol.~30, Springer, Cham, 2019, pp.~197--214. \MR{3931575}

\bibitem{Kressner-Luce-Statti17}
D.~Kressner, R.~Luce, and F.~Statti, \emph{Incremental computation of block
  triangular matrix exponentials with application to option pricing}, Electron.
  Trans. Numer. Anal. \textbf{47} (2017), 57--72. \MR{3707734}

\bibitem{Kurbatov-Kurbatova-LMA16}
V.~G. Kurbatov and I.~V. Kurbatova, \emph{Computation of a function of a matrix
  with close eigenvalues by means of the {N}ewton interpolating polynomial},
  Linear and Multilinear Algebra \textbf{64} (2016), no.~2, 111--122.
  \MR{3434507}

\bibitem{Kurbatov-Kurbatova-CMAM18}
\bysame, \emph{Computation of {G}reen's function of the bounded solutions
  problem}, Comput. Methods Appl. Math. \textbf{18} (2018), no.~4, 673--685.
  \MR{3859261}

\bibitem{Kurbatov-Kurbatova-OM19}
\bysame, \emph{Green's function of the problem of bounded solutions in the case
  of a block triangular coefficient}, Operators and Matrices \textbf{13}
  (2019), no.~4, 981--1001.

\bibitem{Kurbatov-Kurbatova-EMJ20}
\bysame, \emph{An estimate of approximation of a matrix-valued function by an
  interpolation polynomial}, Eurasian Math. J. \textbf{11} (2020), no.~1,
  86--94.

\bibitem{Kurbatov-Kurbatova-Oreshina}
V.~G. Kurbatov, I.~V. Kurbatova, and M.~N. Oreshina, \emph{Analytic functional
  calculus for two operators}, April 2016, arXiv: 1604.07393.

\bibitem{Kurbatov-Oreshina04}
V.~G. Kurbatov and M.~N. Oreshina, \emph{Interconnect macromodelling and
  approximation of matrix exponent}, Analog Integrated Circuits and Signal
  Processing \textbf{40} (2004), no.~1, 5--19.

\bibitem{Kurbatova-Pechkurov20}
I.~V. Kurbatova and A.~V. Pechkurov, \emph{Representations of {G}reen's
  function of the bounded solutions problem for a differential-algebraic
  equation}, Banach J. Math. Anal. \textbf{14} (2020), no.~3, 707--736.
  \MR{4123307}

\bibitem{Mathias93a}
R.~Mathias, \emph{Approximation of matrix-valued functions}, SIAM J. Matrix
  Anal. Appl. \textbf{14} (1993), no.~4, 1061--1063. \MR{1238920}

\bibitem{McCurdy-Ng-Parlett84}
A.~McCurdy, K.~C. Ng, and B.~N. Parlett, \emph{Accurate computation of divided
  differences of the exponential function}, Math. Comp. \textbf{43} (1984),
  no.~168, 501--528.

\bibitem{Moler-Van_Loan78}
C.~Moler and Ch.~F. Van~Loan, \emph{Nineteen dubious ways to compute the
  exponential of a matrix}, SIAM Rev. \textbf{20} (1978), no.~4, 801--836.

\bibitem{Moler-Van_Loan03}
\bysame, \emph{Nineteen dubious ways to compute the exponential of a matrix,
  twenty-five years later}, SIAM Rev. \textbf{45} (2003), no.~1, 3--49
  (electronic).

\bibitem{Parlett74}
B.~N. Parlett, \emph{Computation of functions of triangular matrices}, DTIC
  Document Memorandum ERL-M481, Berkeley, Electronics Research Laboratory,
  College of Engineering, University of California, November 1974.

\bibitem{Parlett76}
\bysame, \emph{A recurrence among the elements of functions of triangular
  matrices}, Linear Algebra and Appl. \textbf{14} (1976), no.~2, 117--121.
  \MR{0448846}

\bibitem{Paterson-Stockmeyer73}
M.~S. Paterson and L.~J. Stockmeyer, \emph{On the number of nonscalar
  multiplications necessary to evaluate polynomials}, SIAM J. Comput.
  \textbf{2} (1973), 60--66. \MR{314238}

\bibitem{Phillips03}
G.~M. Phillips, \emph{Interpolation and approximation by polynomials}, CMS
  Books in Mathematics/Ouvrages de Math\'{e}matiques de la SMC, vol.~14,
  Springer-Verlag, New York, 2003. \MR{1975918}

\bibitem{Sastre18}
J.~Sastre, \emph{Efficient evaluation of matrix polynomials}, Linear Algebra
  Appl. \textbf{539} (2018), 229--250. \MR{3739406}

\bibitem{Sastre-Ibane-Alonso-Peinado-Defez17}
J.~Sastre, J.~Ib\'{a}\~{n}ez, P.~Alonso, J.~Peinado, and E.~Defez, \emph{Two
  algorithms for computing the matrix cosine function}, Appl. Math. Comput.
  \textbf{312} (2017), 66--77. \MR{3660804}

\bibitem{Schmelzer-Trefethen07}
Th. Schmelzer and L.~N. Trefethen, \emph{Computing the gamma function using
  contour integrals and rational approximations}, SIAM J. Numer. Anal.
  \textbf{45} (2007), no.~2, 558--571. \MR{2300287}

\bibitem{Seshu-Reed61:eng}
S.~Seshu and M.~B. Reed, \emph{Linear graphs and electrical networks},
  Addison--Wesley Series in the Engineering Sciences, Addison--Wesley
  Publishing Co., Inc., Reading, Mass.--London, 1961. \MR{0147120}

\bibitem{Varga-Karpenter86}
R.~S. Varga and A.~Dzh. Karpenter, \emph{A conjecture of {S}. {B}ernstein in
  approximation theory}, Mat. Sb. (N.S.) \textbf{129(171)} (1986), no.~4,
  535--548. \MR{842399}

\bibitem{Vlach-Singhal94}
J.~Vlach and K.~Singhal, \emph{Computer methods for circuit analysis and design
  (electrical engineering)}, second ed., Van Nostrand Reinhold
  Electrical/Computer Science and Engineering Series, Kluwer Academic
  Publishers, New York, 1993.

\bibitem{Voevodin-Kuznetsov:rus-eng}
V.~V. Voevodin and Yu.~A. Kuznetsov, \emph{{M}atritsy i vychisleniya
  [{M}atrices and computations]}, ``Nauka'', Moscow, 1984, (in Russian).
  \MR{758446}

\bibitem{Wilkinson88}
J.~H. Wilkinson, \emph{The algebraic eigenvalue problem}, Monographs on
  Numerical Analysis, The Clarendon Press, Oxford University Press, New York,
  1988. \MR{950175}

\bibitem{Wolfram}
S.~Wolfram, \emph{The {M}athematica book}, fifth ed., Wolfram Media, New York,
  2003.

\end{thebibliography}
\end{document}